# Dual Lyapunov-based Synchronization Control of Rössler System*

Alkım Gökçen[1], Savaş Şahin[1], Mahmut Kudeyt[2], Swapnil Tripathi[3] and Özkan Karabacak [3]

[1] Izmir Katip Çelebi University, Izmir 35620, Türkiye
[2] Işık University, Istanbul 34980, Türkiye
[3] Kadir Has University, Istanbul 34083, Türkiye

**Abstract.** This paper proposes a novel approach for the synchronization problem of nonlinear dynamical systems, integrating dual Lyapunov stability analysis with polynomial optimization. A comprehensive review of the relevant scientific literature on synchronization methods is conducted, with a particular focus on classical Lyapunov-based methods for chaotic systems. In this study, the Rössler system is synchronized by employing dual Lyapunov-based closed-loop synchronization method. This method uses semidefinite programming and sum-of-squares polynomials to compute a nonlinear state feedback function which synchronize a chaotic system to a selected reference model. It is aimed that chaotic behavior is destroyed and, instead, a limit cycle becomes attracting. Simulation works are performed for randomly selected 100 different initial conditions to show that synchronization process is successfully performed. Furthermore, bifurcation diagrams and phase portraits are evaluated to analyze the system dynamics. The paper discusses results and how new constraints should be employed and adapted to more complex systems.

## 1 Introduction

Synchronization of the chaotic systems is one of the important ways to practice nonlinear complexities in applications. Chaos control and synchronization of chaotic systems hold an important place is the literature due to their roles in secure communication [1], biology [2], cryptology [3] and finance [4]. Therefore, there are different methods to perform synchronization for chaotic systems including linear and nonlinear feedback design [5,6], adaptive control [7], coupling [8] and online learning control [9]. Lie et al. [10] investigate synchronization algorithms for chaotic systems using a Lyapunov candidate function to introduce a linear feedback controller to achieve synchronization.

* Presented at the International Interdisciplinary Chaos Symposium on Chaos and Complex Systems (SCCS 2025), Istanbul, Türkiye. A version of this work has been accepted for publication in the conference proceedings and will appear in Chaos and Complex Systems: Proceedings of the 6th International Interdisciplinary Chaos Symposium (Springer Cham)

The cornerstone of the Mata-Machuca et al. [11] is to employee Lyapunov stability theory which maintains asymptotic convergence of both synchronization error and the estimation error for uncertain parameters of master/slave systems. Thus, an adaptive synchronization is achieved. Saeed et al. [12] presented a new three-dimensional chaotic system and studied its adaptive synchronization control by constructing a Lyapunov function comprising parameter estimation error and state variables. Herein, they achieved an adaptive control law ensuring synchronization between master and slave systems. Liu et al. [13] utilizes a drive-response synchronization method which rooted in the Lyapunov stability theory. A memristor-based chaotic system is synchronized by constructing a Lyapunov function including synchronization error terms, and method is validated via numerical simulations.

Dual Lyapunov theorem (DLT), introduced by Rantzer [14], characterizes the almost global stability of a dynamical system through a rational density function which increases along the system trajectories. Characterizing almost global stability provides a systematic framework for solving feedback control problems of the systems even beyond the Euclidean domains. DLT has been generalized to time-varying systems [15], switched systems [16], discrete systems [17] and stochastic systems [18]. Leveraging the results from earlier studies, polynomial optimization-based DLT works are implementable for various control engineering problems.

In this study, a nonlinear feedback function based on DLT is computed for synchronization control of the Rössler equation, which is a chaotic system, using sum-of-squares (SOS) polynomials. Initially, complex dynamics of the Rössler system is investigated via phase portrait, bifurcation diagram and Lyapunov exponents to determine parameter conditions. Then, a closed-loop system comprising the master and slave systems is constructed. Herein, the Rössler equation exhibiting chaotic behavior is selected as the slave system, while another Rössler equation exhibiting limit cycle behavior is selected as the master system. Finally, in order to define dual Lyapunov condition, a tracking error-based density function is selected to maintain synchronization between master and slave systems. With the help of the polynomial optimization, a controller is computed and validated through the closed-loop system simulations. The remaining parts of the study are presented under three main sections. The second section explains the theorems and approaches used in this study. Synchronization, evaluation of the DLT together with the synchronization of the Rössler equation is presented in the third section, and simulation results are shared. Lastly, in the fourth section, results of the work and future studies are discussed.

# 2 Mathematical Background: Dual Lyapunov Theorem and Sum-of-Squares Polynomials

The theorems and approaches of this work are explained in this section.

## 2.1 Dual Lyapunov Theorem

A nonlinear dynamical system can be modelled as:

$$\dot{x} = f(x), f: X \to X \tag{1}$$

State variable $x$, defined in $\mathbb{R}^n$, can be analyzed in the perspective of global stability by employing the dual Lyapunov stability criterion [14]. Let $f(x)$ be a continuous and nonlinear system where $f(0) = 0$. If there exist a positive rational function $\rho \in C^1(\mathbb{R}^n - \{0\}, \mathbb{R})$ that is nonnegative and satisfies $[\nabla \cdot (f\rho)(x)] > 0$ (and is integrable on the set of $\{x \in \mathbb{R}^n: |x| > 1\}$), then the system in (1) is almost globally stable. That means, for almost all initial conditions, the solutions of the closed-loop system converge to the zero. Consequently, the DLT can be applied to investigate global stability and to construct nonlinear state feedback functions which can stabilize the equilibrium point.

## 2.2 Sum-of-Squares Polynomials

A polynomial $p$ is called SOS polynomial if it can be written as the sum-of-squares of other polynomials $p_1, p_2, \dots, p_m$, i.e.,

$$p(x) = \sum_{i=1}^{m} p_i^2(x) \tag{2}$$

Representing nonnegative polynomials in SOS form makes it possible to verify nonnegativity of a rational function $\rho$. SOS polynomials are dense in the set of nonnegative polynomials, so approximating nonnegative polynomials using SOS polynomials is the key step. Herein, to solve this step-in practice, the SOS conditions are transformed into linear matrix inequalities (LMI), which can then be solved via semidefinite programming (SDP) methods [19].

# 3 Dual Lyapunov Theorem-based Synchronization Control and Rössler Equation

This section presents a DLT-based nonlinear state feedback function design, Rössler equation and its synchronization study.

## 3.1 Rössler Equation

Rössler equation is proposed by Otto E. Rössler in 1976 [20] as a system which has fixed-point, limit cycle or chaotic behavior under certain parameter conditions. This paper delves into the mathematical structure and properties of the Rössler equation, revealing how change of parameters and initial conditions offers insights into broader phenomena such as synchronization. Rössler system is defined as:

$$\begin{aligned} \dot{x}_1 &= -x_2 - x_3 \\ \dot{x}_2 &= x_1 + a x_2 \\ \dot{x}_3 &= x_1 x_3 + b x_1 - c x_3 \end{aligned} \tag{3}$$

Herein, $a$, $b$ and $c$ are system parameters. It can be seen from (3) that the state equation $\dot{x}_1$ and $\dot{x}_2$ are linear, $\dot{x}_3$ is has one quadratic term $x_1 x_3$ which may cause oscillations in the set-of-solution and generates chaotic complex behaviors. If parameters $a$ and $b$ are fixed ($a = b = 0.2$) and $c$ is selected as control parameter (varies between 0.7 and 8.0), the bifurcation diagram can be analyzed to investigate dynamics of the Rössler system (Fig. 1).

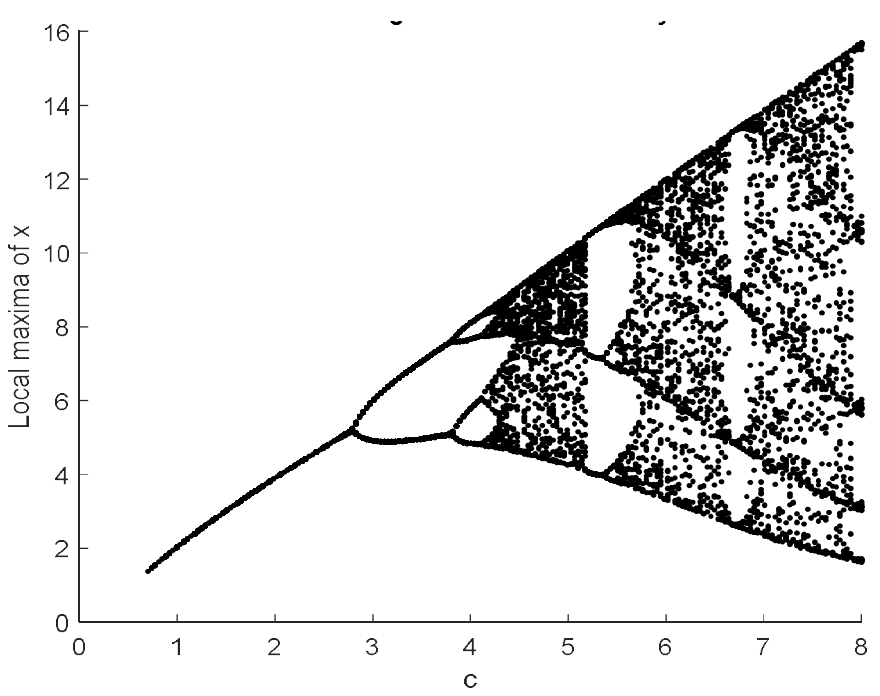

**Fig. 1.** Bifurcation diagram of the Rössler system for $[x_1(0), x_2(0), x_3(0)] = [1,0,0]$.

According to the bifurcation diagram, Rössler system has periodic behavior for $c < 2.8$ and periodicity doubles for larger selections of $c$. For a selection of $c \geq 4.2$ reveals chaotic dynamics in the set-of-solution. Figure 2 represents the phase portrait of the Rössler system for chaotic parameter conditions.

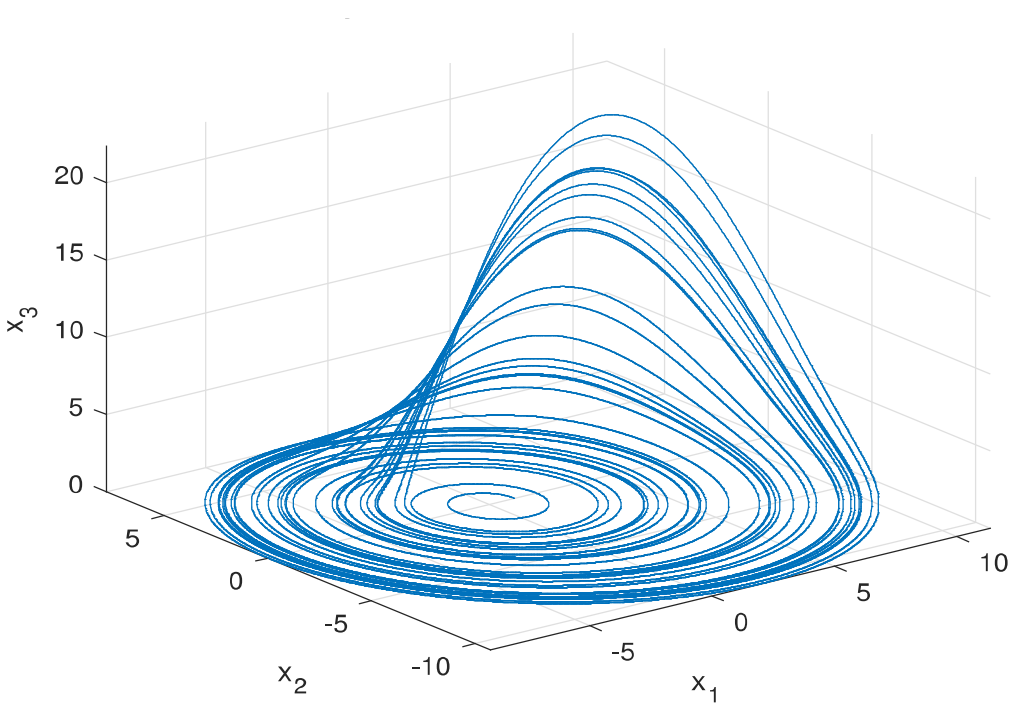

**Fig. 2.** Phase portrait of the Rössler system for $(a, b, c) = (0.2, 0.2, 5.7)$.

### 3.2 Dual Lyapunov-based State Feedback Controller Design and Synchronization

Synchronization process refers to the method aligning the states and/or outputs of two or more dynamics systems over time through appropriate coupling or feedback function. It plays significant role in nonlinear and chaotic systems. By designing a suitable control input, generally applied to slave system, synchronization achieves that its dynamics closely follows the master system dynamics. For this purpose, a closed-loop system comprising master and slave systems can be written as follows:

$$\begin{bmatrix} \dot{x}_1 \\ \dot{x}_2 \\ \dot{x}_3 \\ \dot{x}_{d_1} \\ \dot{x}_{d_2} \\ \dot{x}_{d_3} \end{bmatrix} = \begin{bmatrix} -x_2 - x_3 + u \\ x_1 + 0.2x_2 \\ x_1x_3 + 0.2x_1 - 5.7x_3 \\ -x_{d_2} - x_{d_3} \\ x_{d_1} + 0.2x_{d_2} \\ x_{d_1}x_{d_3} + 0.2x_{d_1} - 2.5x_{d_3} \end{bmatrix} \tag{4}$$

The designed closed-loop system in (4) consists of a master system $x_d$ in limit cycle which represents the desired dynamics, and a slave system $x$ in chaos. The control input $u$ is designed as a state feedback function considering the dual Lyapunov positivity constraints in MATLAB environment. Herein, rational density function is selected in a form which characterizes the almost global stability of the error system defined by $\dot{e} = \dot{x} - \dot{x}_d$:

$$\rho = \frac{1}{e_1^2 + e_2^2 + e_3^2} \tag{5}$$

Thus, satisfying the DLT condition maintains convergence to zero with time for denominator of the (5) and ensures synchronization. YALMIP and SOSTOOLS libraries are employed define constraints in LMI form and SeDuMi SDP solver is used to compute a second order polynomial state feedback function, $u$:

$$u(x) = x^T \begin{bmatrix} 1.120 & 0.565 & -0.020 & -0.300 & -0.565 & -0.555 \\ 0.565 & 0.010 & -0.075 & 0.255 & -0.010 & -0.005 \\ -0.02 & -0.075 & -0.160 & 0.340 & 0.075 & 0.080 \\ -0.300 & 0.255 & 0.340 & -0.520 & -0.255 & 0.235 \\ -0.565 & -0.100 & 0.075 & -0.255 & 0.010 & 0.005 \\ -0.555 & -0.005 & 0.080 & 0.235 & 0.005 & 0.210 \end{bmatrix} x \quad (6)$$

The closed-loop system is simulated in Simulink/MATLAB environment with the computed state feedback function in (6). Herein, 100 randomly selected initial conditions are employed in simulation works, and phase portrait of the slave system is visualized. According to phase portrait (Fig. 3), trajectories of the slave system evolve unison with the master system overtime, and synchronization is achieved.

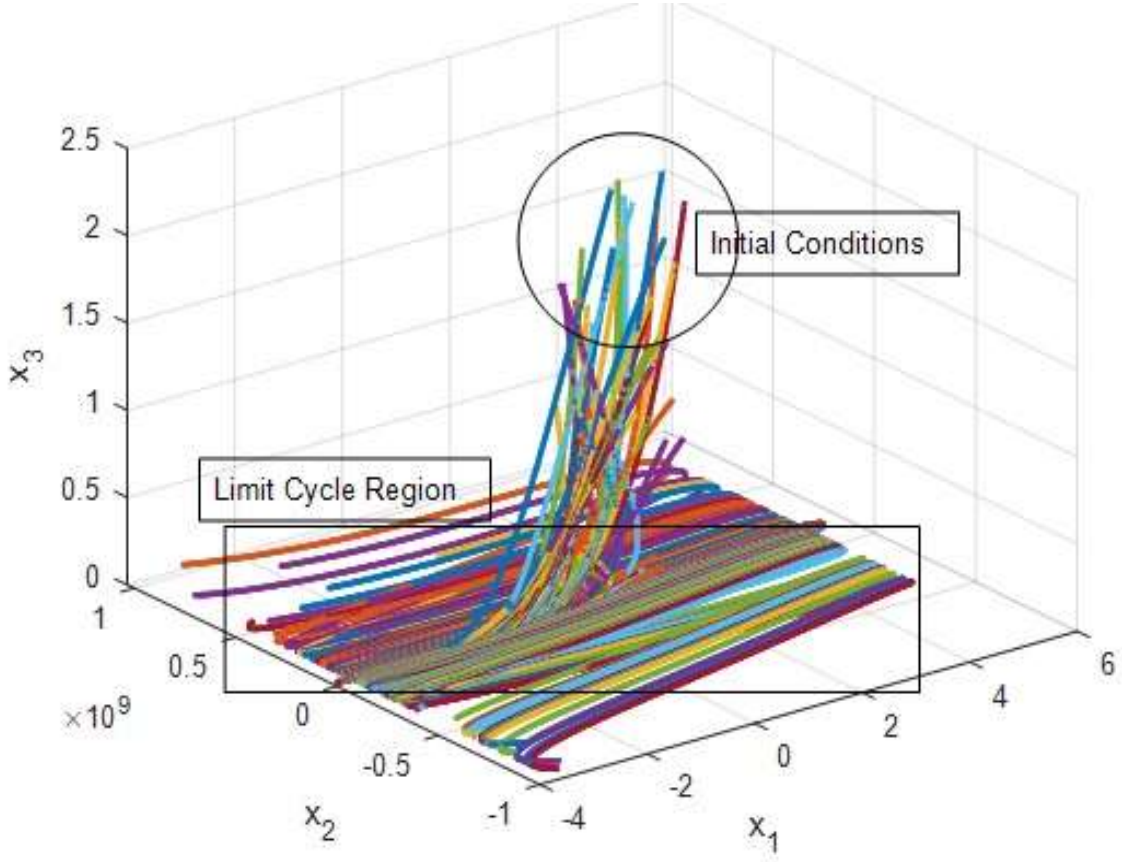

**Fig. 3.** Phase portrait evaluation of the slave Rössler system.

## 4 Conclusion

This study presents a chaotic system synchronization framework for nonlinear dynamical systems using a dual Lyapunov-based method combined with polynomial optimization and semidefinite programming techniques. By leveraging a density function-based stability analysis, the proposed approach enables almost global stability of the closed-loop system. The Rössler equation is employed as a benchmark, and the proposed method successfully derived a state feedback function to ensure synchronization across variously different initial conditions. SeDuMi SDP solver is used to compute the state feedback function, and simulation studies validated the synchronization process. Moving forward, the proposed method can be extended to more intricate systems and refined to incorporate multi-constrained control design.